\documentclass[a4paper,12pt]{amsart}
\usepackage{amscd,amsmath,amssymb,amsxtra,boldsect}
\usepackage{hyperref}
\usepackage[xdvi]{graphicx}

\newcommand{\Q}{\mathbb Q}

\newcommand{\Z}{\mathbb Z}
\newcommand{\Ga}{\alpha}
\newcommand{\Gb}{\beta}
\newcommand{\Gs}{\sigma}
\newcommand{\p}{\partial}
\newcommand{\rnk}{\operatorname{rk}}
\newcommand{\cal}{\mathcal}
\newcommand{\bfit}{\bfseries\itshape}

\title{Remarks on definition of Khovanov homology}

\author{Oleg Viro}
\dedicatory{Uppsala University, Uppsala, Sweden\break
POMI, St.\ Petersburg, Russia}
\address{Department of Mathematics, Uppsala University,
Box 480, S-751 06 Uppsala, Sweden
}
\email{oleg@math.uu.se}
\begin{document} 
\begin{abstract}

Mikhail Khovanov defined, for a diagram of an oriented classical 
link, a collection of groups numerated by pairs of
integers.   These groups were constructed as homology groups of  
certain chain complexes. The Euler characteristics of these complexes 
are coefficients of the Jones polynomial of the link. The goal of this 
note is to rewrite this construction in terms more friendly to 
topologists. A version of Khovanov homology for framed links is
introduced. For framed links whose Kauffman brackets are involved in a
skein relation, these homology groups are related by an exact sequence.
\end{abstract}
\maketitle

\section{Introduction}\label{s0}

For a diagram $D$ of an oriented link $L$,
Mikhail Khovanov \cite{Kh} constructed a collection of groups
$\cal H^{i,j}(D)$ such that
$$K(L)(q)=\sum_{i,j}q^{j}(-1)^i\dim_{\Q}(\cal H^{i,j}\otimes\Q),$$
where $K(L)$ is a version of the Jones polynomial of $L$.
These groups are constructed as homology groups of chain complexes. 
The primary goal of this note is to rewrite this construction in terms more
pleasant for topologists. 

To some extent this has been done recently by Dror Bar-Natan in his 
preprint \cite{BN}. This came together with progress in
calculation and understanding of Khovanov homology, see \cite{Kh2} and
\cite{Lee}.  The outbreak of activity that happened two years after
the release of  initial paper has proved that stripping Khovanov's
construction of its fancy  formal decorations was useful. 

I hope that a further chewing\footnote{From Bar-Natan's paper: ``Not being 
able to really digest it (construction of Khovanov's invariant -- O.V.) 
we decided to just chew some, and then provide our output 
as a note containing a description of his construction, complete and 
consistent and accompanied by computer code and examples but stripped 
of all philosophy and of all the linguistic gymnastics that is necessary 
for the philosophy but isn't necessary for the mere purpose of having 
a working construction. Such a note may be more accessible than the 
original papers. It may lead more people to read Khovanov at the source, 
and maybe somebody reading such a note will figure out what the Khovanov 
invariants really are. Congratulations! You are reading this note right
now.''} of Khovanov's construction that is presented below  occur to
be useful.  I am grateful to Khovanov,
who made these my  considerations possible by not only the very
initiating of the subject, but also by keeping  (at least
during the Fall of 1998) his text secret, while  giving talks about his
work at various seminars. The rumors reached me, and I tried to figure
out how homology with the Euler characteristic equal to the Kauffman
bracket could be defined. When the preprint became available, I found
that Khovanov's construction can be reformulated in the way that I had
guessed. I could upload then most of the text presented below, but 
hesitated, since I had no real new results based on my chewing. I still
have no them, but feel that a text showing how a transition from Kauffman 
bracket to Khovanov homology could be motivated for a topologist may be
appreciated. Some irresponsible speculations about possible
generalizations are included at the end of Section \ref{s1}.

The Khovanov homology is closer to Kauffman bracket, which is an invariant
of non-oriented, but framed links, rather than to the Jones polynomial, which
is an invariant of oriented, but non-framed links. The corresponding
modification of Khovanov homology is presented in Section
\ref{s2}. This allows us to write down a categorification of the Kauffman
skein relation for Kauffman bracket. The skein relation gives rise to a 
homology sequence.

\section{From Kauffman bracket to Khovanov homology}\label{s1}

Since the Euler characteristic does not change under passing from a
chain complex to its homology, it is natural to expect that $\cal
H^{i,j}$ appear as the homology groups of a chain complex $\cal
C^{i,j}$ such that its polynomial 
Euler characteristic $\sum_{i,j}q^j(-1)^i\rnk \cal C^{i,j}$ is the
Jones polynomial $K(L)(q)$.

Khovanov constructs such a complex starting with the Kauffman state sum
presentation of the Jones polynomial. However, the construction
proceeds as a chain of algebraic constructions of auxiliary objects.
This hides a simple geometric meaning. A natural generators of
Khovanov's complex can be represented as enhancements of the states
from the Kauffman state sum presentation of the Jones polynomial.

\subsection{Kauffman's states of a link diagram}\label{s1.1}
\begin{figure}[htb]
\centerline{\includegraphics{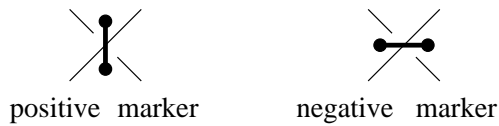}}
\caption{Markers comprising a Kauffman state.}
\label{f-markers}
\end{figure}

Recall that 
for a state sum
representation of the Jones polynomial, Kauffman
\cite{Ka} introduced   the following states of a link
diagram. Each state is a collection of markers. At each crossing of the
diagram there should be a marker specifying a pair of vertical angles.
See Figure \ref{f-markers}.

\subsection{Numerical characteristics of states and
diagrams}\label{s1.2} 
\begin{figure}[htb]
\centerline{\includegraphics{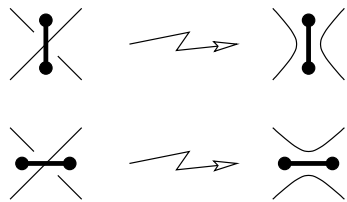}}
\caption{Smoothing of a diagram according to markers.}
\label{f-smoothings}
\end{figure}

Each marker of a Kauffman state has a {\em sign},
which is $+$ (or, rather, $+1$) if the direction of rotation of the 
upper string
towards the lower one through the specified pair of vertical angles is
counter-clock-wise, and $-$ (or $-1$) otherwise. For a state $s$ of a
diagram $D$ denote by $\Gs(s)$ the difference between the numbers of
positive and negative markers. A state of a diagram defines a {\em
smoothing\/} of the diagram: at each of its double points the marked
angles are united in a connected area, see Figure \ref{f-smoothings}.

The result of the smoothing is a collection of circles embedded into the
plane. Denote the union of these circles by $D_s$. The number of the
circles is denoted by $|s|$.

With a crossing point of a diagram of an oriented
link we associate {\em local writhe number} equal to\begin{itemize}
\item $+1$, if at the
point the diagram looks like
$\vcenter{\hbox{\includegraphics{fig/pcross.eps}}}$, and
\item $-1$, if it looks like 
$\vcenter{\hbox{\includegraphics{fig/ncross.eps}}}$.\end{itemize}
The sum of local writhe numbers over all crossing points of a link
diagram $D$ is called the {\it writhe number\/} of $D$ and denoted by $w(D)$.

To a Kauffman state $s$ of an oriented link diagram $D$ we assign a
polynomial
$$V_s(A)=A^{\Gs(s)}(-A)^{-3w(D)}(-A^2-A^{-2})^{|s|},$$
where $w(D)$ is the writhe number of $D$. 

\subsection{Jones Polynomial}\label{s1.3} The sum
of $V_s(A)$ over all states $s$ of the diagram is a version the Jones
polynomial of $L$. It is denoted by $V_L(A)$. 

Khovanov \cite{Kh} uses another variable $q$ instead of $A$. It is 
related to $A$ by $q=-A^{-2}$. Following Khovanov, we denote the 
corresponding version of the Jones polynomial by $K(L)$. It is defined by
$K(L)(-A^{-2})=V_L(A)$.

Certainly, $K(L)$ can be presented as the sum of the corresponding
versions of $V_s(A)$ over all Kauffman states $s$. The summand of
$K(L)$ corresponding to $s$ can be defined via $K(s)(-A^2)=V_s(A)$ or,
more directly, by 
$$K(s)(q)=(-1)^{\frac{w(D)-\Gs(s)}2}q^{\frac{3w(D)-\Gs(s)}2}(q+q^{-1})^{|s|}.$$ 

\subsection{Enhanced Kauffman states}\label{s1.4}
By an {\em enhanced Kauffman state\/} $S$ of an oriented link  diagram
$D$ we shall mean a collection of markers comprising a usual Kauffman
state $s$ of $D$ enhanced by an assignment of a plus or minus sign to
each of the circles obtained by the smoothing of $D$ according to the
markers.

Denote by $\tau(S)$ the difference between the numbers of pluses and
minuses assigned to the circles of $D_s$. Put
$$j(S)=-\frac{\Gs(s)+2\tau(S)-3w(D)}2.$$
Observe that both $\Gs(s)$ and $w(D)$ are congruent modulo 2 to the
number of crossing points. Therefore $j(S)$ is an integer.

\subsection{Monomial state sum representation of Jones
polynomial}\label{s1.5} Recall that
$V_s(A)=A^{\Gs(s)}(-A)^{-3w(D)}(-A^2-A^{-2})^{|s|}$. We can associate
to each component of $D_s$ one of $|s|$ factors $(-A^2-A^{-2})$ of
the right hand side of this formula.

Passing to the enhanced Kauffman states, we associate the summand
$-A^2$ of the sum $-A^2-A^{-2}$ to a component of $D_{|s|}$ equipped
with $+$. To the same component, but equipped with $-$, we associate
the summand $-A^{-2}$.

Thus an enhancement of a Kauffman state gives rise to a choice of a
summand from each of the binomial factors in the formula defining
$V_s(A)$. Now by opening the brackets we obtain the following
presentation for $V_s(A)$:
\begin{multline}\label{SpltVs(A)}
V_s(A)= \\
\sum_{S\in s}A^{\Gs(s)}(-A^2)^{\tau(S)}(-A)^{-3w(D)}= \\
\sum_{S\in s}(-1)^{\tau(S)+w(D)}A^{-2j(S)}.
\end{multline}
Denote by $V_S(A)$ the summand $(-1)^{\tau(S)+w(D)}A^{-2j(S)}$
which corresponds to $S$.
Hence $V_s(A)=\sum_{S\in s}V_S(A)$ and
\begin{multline}\label{V_L(A)}
V_L(A)=\sum_{\substack{\text{Kauffman}\\
\text{ states }s \text{ of }D}} V_s(A)= \\
\sum_{\substack{\text{enhanced Kauffman}\\
\text{ states }S \text{ of }D}} V_S(A)=
\sum_{\substack{\text{enhanced Kauffman}\\
\text{ states }S \text{ of }D}}
(-1)^{\tau(S)+w(D)}A^{-2j(S)}.\end{multline}

\subsection{Change of variable}\label{s1.6}
Let us switch to variable $q$ used by Khovanov. 
\begin{multline} V_S(A)=
(-1)^{\tau(S)+w(D)}A^{-2j(S)}=\\
(-1)^{\tau(S)+w(D)}(A^{-2})^{j(S)}=\\
(-1)^{\tau(S)+w(D)-j(S)}(-A^{-2})^{j(S)}=\\
(-1)^{\tau(S)+w(D)+\frac12\Gs(s)+\tau(S)-\frac32w(D)}q^{j(S)}=\\
(-1)^{\frac{\Gs(s)-w(D)}2}q^{j(S)}=\\
(-1)^{\frac{w(D)-\Gs(s)}2}q^{j(S)}.
\end{multline}
Denote $(-1)^{\frac{w(D)-\Gs(s)}2}q^{j(S)}$ by $K(S)(q)$.
Clearly,
\begin{equation}\label{V_L(q)}
K(L)(q)=
\sum_{\substack{\text{enhanced Kauffman}\\
\text{ states }S \text{ of }D}}
(-1)^{\frac{w(D)-\Gs(s)}2}q^{j(S)}
=\sum_{\substack{\text{enhanced Kauffman}\\
\text{ states }S \text{ of }D}} K(S)(q).\end{equation}

\subsection{Khovanov chain groups}\label{s1.7}
Denote the free abelian group generated by enhanced Kauffman states of
a link diagram $D$ by $\cal C(D)$. Denote by $\cal C^j(D)$ the subgroup of $\cal C(D)$
generated by enhanced Kauffman states $S$ of $D$ with $j(S)=j$. Thus
$\cal C(D)$ is a $\Z$-graded free abelian group:
$$\cal C(D)=\bigoplus_{j\in\Z}\cal C^j(D).$$

For an enhanced Kauffman state $S$ belonging to a Kauffman state $s$ of
a link diagram $D$, put
$$i(S)=\frac{w(D)-\Gs(s)}2.$$
Denote by $\cal C^{i,j}(D)$ the subgroup of $\cal C^j(D)$ generated by enhanced
Kauffman states $S$ with $i(S)=i$.

Notice that as follows from \eqref{V_L(q)},
\begin{equation}\label{chiHij}
K(L)(q)=\sum_{j=-\infty}^{\infty}q^j\sum_{i=-\infty}^{\infty}
(-1)^i\rnk \cal C^{i,j}(D).
\end{equation}

{\bfit Remark.\/}
The Kauffman choice of variable ($A$ instead of $q$) in the Jones 
polynomial discussed above suggests another choice of the second 
grading of $\cal C(D)$. In
fact, the choices incorporate already two different $\Z_2$-grading:
$\tau(S)+w(D)\pmod 2$ and $\frac{w(D)-\Gs(s)}2\pmod 2$. However
these $\Z_2$-grading differ just by the reduction modulo 2 of the first
one, $$\tau(S)+w(D)-\frac{w(D)-\Gs(s)}2=j(S)\mod2.$$ 
One can prove that the Kauffman choice of variable gives rise to the
same homology groups.

\subsection{Differential}\label{s1.8}
Let $D$ be a diagram of an oriented link $L$. In \cite{Kh} Khovanov
defined in $\cal C^{i,j}(D)$ a differential of bidegree $(1,0)$ and proved
that its homology groups $\cal H^{i,j}(D)$ do not depend on $D$. The
construction of the differential depends on ordering of crossing points
of $D$, although the homology groups do not. Suppose that the crossing
points of $D$ are numerated by natural numbers $1,\dots,n$.

Khovanov's description of the differential is somewhat complicated by 
several auxiliary algebraic constructions. I give here its simplified 
version. I just describe the matrix elements. In the context of chain 
complex the matrix elements are traditionally called {\em incidence 
numbers}. For enhanced Kauffman states $S_1$ and $S_2$, denote their 
incidence number by $(S_1:S_2)$.

Recall that $\cal C^{i,j}(D)$ is generated by enhanced Kauffman states. 
Thus an incidence number is a function of two enhanced Kauffman 
states, say $S_1$ and $S_2$ (which are generators of $\cal C^{i,j}(D)$ and 
$\cal C^{i+1,j}(D)$, respectively). Enhanced Kauffman states with a nonzero 
incidence number are said to be {\it incident} to each other. Pairs of 
incident states satisfy natural restrictions. Surprisingly, these 
restrictions give exact description of the set of incident 
states: each pair of enhanced Kauffman states which is not thrown away
by the restrictions consists of incident states.

The first restriction emerges from our desire to have a differential 
of bidegree $(1,0)$. Since the differential preserves $j$ and
increases $i$ by one, $j(S_1)=j(S_2)$ and $i(S_2)=i(S_1)+1$. Recall
that $i(S)=\frac{w(D)-\Gs(s)}2$. Therefore $i(S_2)=i(S_1)+1$ implies
that $\Gs(s_2)=\Gs(s_1)-2$, in other words the number of negative
markers of $S_2$ is greater by one than the number of negative markers
of $S_1$.

It is natural to enforce this numerological restriction in the 
following way:

{\it The incidence number is zero, unless only at one crossing point
of $D$ the markers of $S_1$ and $S_2$ differ and at this crossing the
marker of $S_1$ is positive, while the marker of $S_2$ is negative.}

In our description of the incidence number of $S_1$ and $S_2$ we
will assume that this is the case. The crossing point where the 
markers differ is called the {\em difference point} of $S_1$ and 
$S_2$. Let it have number $k$. 

Since exactly at one crossing the markers differ, $D_{S_2}$
is obtained from $D_{S_1}$ by a single oriented Morse
modification of index 1. Hence $|S_1|-|S_2|=\pm1$. In other words
either $D_{S_2}$ is obtained from $D_{S_1}$ by joining two circles or
by splitting a circle of $D_{S_1}$ into two circles.

Here is the next natural restriction on incident states:

{\em The incidence number of $S_1$ and $S_2$ vanishes, unless the
common circles of $D_{S_1}$ and $D_{S_2}$ have the same signs.}

This, together with equality $j(S_2)=j(S_1)$, gives a strong 
restriction
also on the signs of the circles of $D_{S_l}$ adjacent to the $k$th
crossing point. Indeed, 
\begin{multline*}
j(S_2)=\frac{3w(D)-\Gs(s_2)-2\tau(S_2)}2\\=
j(S_1)=\frac{3w(D)-\Gs(s_1)-2\tau(S_1)}2=\frac{3w(D)-\Gs(s_2)-2-2\tau(S_1
)}2,\end{multline*} 
thus $\tau(S_2)=\tau(S_1)+1$. 

\begin{figure}[thbp]
\centerline{\includegraphics{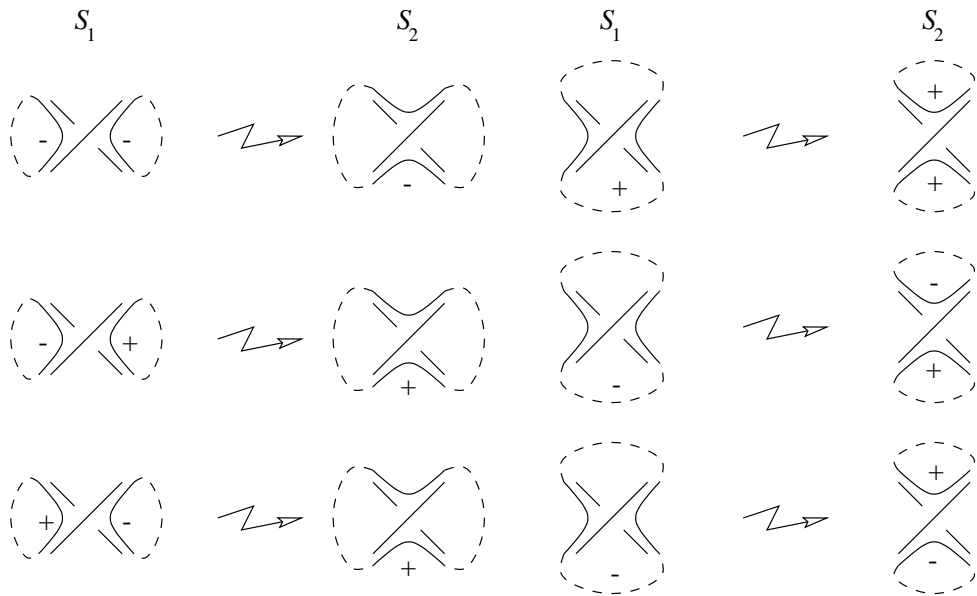}}
\caption{Pairs of incident enhanced Kauffman states.
Dotted arcs show how the fragments of $D_s$ at a crossing point are
connected in the whole $D_s$.}
\label{f-Xinc}
\end{figure}

Now we can easily list all the
situations which satisfy these restrictions (see Figure \ref{f-Xinc}):
\begin{enumerate}
\item If $|S_2|=|S_1|-1$ and
both joining circles of $D_{S_1}$ are negative then the resulting
circle of $S_2$ should be negative.
\item If $|S_2|=|S_1|-1$ and
the joining circles of $D_{S_1}$ have different signs then the resulting
circle of $S_2$ should be positive.
\item If the joining circles of $D_{S_1}$
are positive, none $S_2$ can be incident.
\item If $|S_2|=|S_1|+1$ and the
splitting circle of $D_{S_1}$ is positive then both of the circles of
$D_{S_2}$ obtained from it should be positive.
\item If $|S_2|=|S_1|+1$ and the
splitting circle of $D_{S_1}$ is negative then the circles of
$D_{S_2}$ obtained from it should be of different signs.
\end{enumerate}

In each of the cases listed above the incidence number $(S_1:S_2)$ is 
$\pm1$. The sign depends on the ordering of the crossing points: it is 
equal to $(-1)^t$ where $t$ is the number of negative markers in $S_1$ 
numerated with numbers greater than $k$. 

This sign is needed to make the differential satisfying identity
$d^2=0$. However, the proof of this identity requires some routine check.
At first glance, Khovanov \cite{Kh} escaped it. In fact, it is hidden in
the checking (which was left to the reader)  that $F$ is a functor, see
\cite{Kh}, Section 2.2.

\subsection{Enhanced Kauffman states with polynomial
coefficients}\label{s1.9} Khovanov constructed not only groups
$\mathcal H^{i,j}$,  but also graded modules $H^i$ over the ring
$\Z[c]$ of polynomials with integer coefficients in variable $c$ of
degree 2. The grading is a representation of $H^i$ as a direct sum of
abelian subgroups $H^{i,j}$ such that multiplication by $c$ in $H^i$
gives rise to a homomorphism $H^{i,j}\to H^{i,j+2}$.

To construct homology groups $H^{i,j}$, let us define the corresponding
complex of graded $\Z[c]$-modules $C^i$. The module $C^i$ is the sum of
its subgroups $C^{i,j}$. The group $C^{i,j}$ is generated by  formal
products $c^kS$, where $k\ge0$ and $S$ is an enhanced Kauffman states
$S$ with  $i(S)=i$ and $j(S)=j-2k$. 

\begin{figure}[htb]
\centerline{\includegraphics{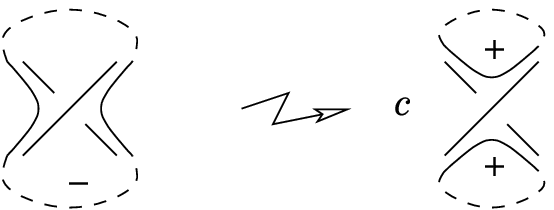}}
\caption{}
\label{c-inc}
\end{figure}

The differential is defined almost in the same way. The states which
were adjacent above are adjacent here, as well as their products by the
same power of $c$. Products of enhanced Kauffman states by different 
powers of $c$ are not adjacent, besides the following  situations: 
$c^{k+1}S_2$ is adjacent to $c^kS_1$, if the markers of $S_1$ and $S_2$
differ exactly at one point, where the marker of $S_1$ is positive and
the marker of $S_2$ negative, the signs of $S_1$ and $S_2$ on the common
circles of $D_{S_1}$ and $D_{S_2}$ are the same, $|S_2|=|S_1|+1$,
the splitting circle of $D_{S_1}$ is negative, and the circles of
$D_{S_2}$ obtained from it are positive, see Figure \ref{c-inc}, where
this situation is shown symbolically in the style of Figure \ref{f-Xinc}.

Each group $C^{i,j}$ is finitely generated, but there are infinitely
many non-trivial groups. Groups $H^{i,j}$ with fixed $i$ and
sufficiently large $j$ are isomorphic to each other. 

\subsection{Speculations}\label{s1.10}
How can one understand construction of Khovanov homology? By
understanding I mean a formulation which would give, true or false,
indication to possible generalizations. The general program of
categorification pictured by Khovanov in \cite{Kh} is an attractive
approach to the problem of constructing combinatorial counterparts for
Donaldson and Seiberg-Witten invariants, although it does not promise a
fast and technically easy track.

In the construction above the specific of Kauffman bracket is used so
much that it is difficult to imagine a geralization to other polynomial
link invariants. A state sum over smoothings along Kauffman markers does not
appear in formulas representing any other quantum invariant. Of course,
there are many state sum formulas. One can try to decompose any of them
to a sum of monomials and consider decompositions as chains. Attach to them
appropriate dimensions and define differentials such that the homology
groups would be invariant under Reidemeister moves. Success of the
whole project  depends on quality of several guesses. The most difficult 
of the guesses seems to be the choice of differentials. It may involve and 
rely on a delicate analysis of the effects of Reidemeister moves. However 
a wrong choice at the previous step can make any efforts at this last step
unfruitful.

An additional hint comes from Khovanov's mention of
Lusztig's canonical bases. Indeed, signs which are attached to circles 
in a smoothened link diagram can be replaced with orientations and an
orientation was interpreted by Frenkel and Khovanov \cite{FK} as an
element of Lusztig's canonical basis. Smoothings at crossing points
correspond then to monomials of the entries of the $R$-matrix. The 
decomposition of the entries of $R$-matrix to monomials not as clear as
it seams: in the only known case considered above a zero entry of 
$R$-matrix is decomposed to two non-zero monomials canceling each
other. However, the main guess is still to be done after the entries of
the $R$-matrix in crystal bases are decomposed to monomials: one should
find a differential that would give homology invariant under
Reidemeister moves. 

If the guesses made above about the r\^oles of crystal bases and
monomials of $R$-matrices are correct, $1+1$-dimensional TQFTs would not
appear in the future more involved categorifications. Instead of cobordisms
of 1-manifolds, there appear probably cobordisms of 4-valent 
graphs with some additional structure. 

How can one estimate chances for further categorification of quantum
topology? Still, they do not seem to be clear, despite of high optimism
of Khovanov.As for chain groups, the observations made above
are encouraging. However, we are not even able to discuss the choice 
of differentials on a comparable level of generality.

It would be very interesting, especially for prospective 4-dimensional
applications, to find a categorification of face state models.

\section{Frame version}\label{s2}

\subsection{Kauffman bracket versus Jones polynomial}\label{s2.1}
Involvement of an orientation of link in the definition of the Jones
polynomial in Section \ref{s1.3} is easy to localize. It is limited to
to the writhe number $w(D)$ of $D$. If we remove $w(D)$, we get the
{\it Kauffman bracket\/} of diagram $D$. To eliminate fractional powers,
we stay with Kauffman's variable $A$. Thus Kauffman bracket of a link
diagram $D$ is 
$$\langle D\rangle=
\sum_{\substack{\text{ Kauffman }\\ \text{ states }
s\text{ of }D}}A^{\Gs(s)}(-A^2-A^{-2})^{|s|}.$$
It is not invariant under first Reidemeister moves. Hence it is not 
an ambient isotopy invariant of link. However it is invariant with respect to 
second and third Reidemeister moves. 

A link diagram defines the blackboard framing on the link. The
isotopies which accompany second and third Reidemeister moves can be
extended to the  isotopy of the corresponding framings. Moreover,
embedding of the diagrams to the 2-sphere makes this relation
two-sided: link diagrams on $S^2$ can be converted to each other by 
a sequence of second and third Reidemeister moves iff the corresponding
links with blackboard framings are isotopic in the class of framed
links. Since the Kauffman bracket can be defined for a link diagram on
2-sphere, this makes the Kauffman bracket invariant for non-oriented
framed links.

The Kauffman bracket is categorified below.

\subsection{Framed Khovanov homology}\label{s2.2}
For an enhanced Kauffman state of a link diagram $D$, put 
$$J(S)=\Gs(s)+2\tau(S)\quad\text{ and }\quad I(S)=\Gs(s). $$
If this is an oriented link, and hence $w(D)$ is defined, then
$$J(S)=3w(D)-2j(S)\quad\text{ and }\quad I(S)=w(D)-2i(S).$$ 
In terms of $J(S)$ and $I(S)$, the Kauffman bracket is expressed as
follows: 
$$\langle D\rangle=\sum_{\substack{\text{enhanced Kauffamn}\\
\text{states } S \text{ of }  D}}(-1)^{\frac{I(S)}2}A^{J(S)}.$$

Denote the free abelian group generated by enhanced Kauffman states $S$ of
$D$ with $I(S)=i$ and $J(S)=j$ by $\cal C_{i,j}(D)$. If one orients the
link, the Khovanov chain groups $\cal C^{i,j}(D)$ appear, and
$$\cal C_{i,j}(D)=\cal C^{\frac{w(D)-i}2,\frac{3w(D)-j}2}(D).$$
Under this identification, the differentials of the Khovanov complex
turn into differentials
$$\p:\cal C_{i,j}(D)\to \cal C_{i-2,j}(D)$$
(the construction of the differentials does not involve orientation of
the link, hence it does not matter what orientation is used). Denote the
homology group of the complex obtained by $H_{i,j}(D)$.

\subsection{Skein homology sequence}\label{s2.3} Let $D$ be a link
diagram, $c$ its crossing and $D_+$, $D_-$ be link diagrams obtained
from $D$ by smoothing at $c$ along positive and negative markers,
respectively. As is well-known, the Kauffman brackets of $D$, $D_+$ and
$D_-$ are  related as follows
$$\langle D\rangle=A\langle D_+\rangle+A^{-1}\langle D_-\rangle.$$
This equality is called the {\it Kauffman skein relation.\/} It allows one
to define the Kauffman bracket up to a normalization (which can be
done by fixing the Kauffman bracket of the unknot).

Let us categorify this skein relation.
Consider map 
$$\Ga:\cal C_{i,j}(D_-)\to \cal C_{i-1,j-1}(D)$$ 
which sends an enhanced Kauffman state $S$ of $D_-$ to the enhanced Kauffman
state of $D$ smoothing along which coincides with the smoothing of
$D_-$ along $S$ and signs of the ovals are the same, too.
The collection of these maps is a homomorphism of complex, that is they
commute with $\p$, provided $c$ is the last crossing of $D$ in the 
ordering of crossings which is used in the construction of $\p$.
Indeed, the incidence coefficients are the same for an enhanced
Kauffman states of $D_-$ and its images in $D$, and the latter cannot
contain in the boundary a state with positive marker at $c$.

Now consider map 
$$\Gb:\cal C_{i,j}(D)\to \cal C_{i-1,j-1}(D_+)$$
which sends each enhanced Kauffman state with negative marker at $c$ to
0 and each enhanced Kauffman state with positive marker at $c$ to the
enhanced Kauffman state of $D_+$ with the same smoothing and signs of
the ovals. This is again a complex homomorphism.

Homomorphisms $\Ga$ and $\Gb$ form a short exact sequence of complexes:
$$\begin{CD}
0@>>>\cal C_{*,*}(D_-)@>{\Ga}>> \cal C_{*-1,*-1}(D)@>{\Gb}>> 
\cal C_{*-2,*-2}(D_+) @>>> 0
\end{CD}
$$
It induces collection of long homology sequences:
$$
\begin{CD}
@>{\p}>>H_{i,j}(D_-)@>{\Ga_*}>> H_{i-1,j-1}(D)@>{\Gb_*}>> 
H_{i-2,j-2}(D_+)\\ @>{\p}>>H_{i-2,j}(D_-)
@>{\Ga_*}>> H_{i-3,j-1}(D)@>{\Gb_*}>>H_{i-4,j-2}(D_+) @>{\p}>>
\end{CD}
$$
A special case of this sequence, which relates the groups of connected
sum and disjoint sum of knots, can be found in Section 7.4 of
Khovanov's paper \cite{Kh}. This special case is the only one which can
be formulated for the original version of homology depending on
orientations of links.

I was not able to categorify the other skein relation, the one which involves
the Jones polynomial of {\it oriented\/} links. It follows from a couple of
skein relations of the type considered above.


\begin{thebibliography}{9999999}

\bibitem{BN} Dror Bar-Natan,  Khovanov's categorification of the Jones
polynomial, arXiv: math.QA /0201043.

\bibitem{FK} I.~B.~Frenkel, M.~Khovanov, Canonical bases in tensor
products and graphical calculus for $U_q(\mathfrak{sl}_2)$, {\it Duke
Math. Journal,\/} v. 87, (1997), 409--480. 
\bibitem{Ka}  L.H.Kauffman, State models and the Jones polynomial,
{\it Topology,\/} {\bf 26:3} (1987), 395--407.

\bibitem{Kh} Mikhail Khovanov, A categorification of the Jones
polynomial, {\it Duke Math. J.,\/} 101 (3), 359--426, 1999,
arXiv: math.QA /9908171.

\bibitem{Kh1} Mikhail Khovanov, A functor-valued invariant of tangles,
arXiv: math.QA /0103190

\bibitem{Kh2} Mikhail Khovanov, Patterns in knot cohomology I,
arXiv: math.QA /0201306.

\bibitem{Lee} E.~S.~Lee, The support of the Khovanov's invariants for
alternating knots, arXiv: math.GT /0201105.

\end{thebibliography}
\end{document}